\documentclass{article}
\usepackage{graphicx} 
\usepackage{amsfonts,amssymb,amsmath,amsthm}
\usepackage{xcolor}
\usepackage{appendix}
\usepackage{nicefrac}
\usepackage{placeins}
\usepackage{bbm}
\usepackage{physics}
\usepackage[normalem]{ulem}
\usepackage{tikz, pgfplots}
\usepackage{url}

\numberwithin{equation}{section}
\allowdisplaybreaks

\title{Data-driven MPC with stability guarantees using extended dynamic mode decomposition\thanks{K.~Worthmann gratefully acknowledges funding by the Deutsche Forschungsgemeinschaft (DFG, German Research Foundation) -- {Project-ID 507037103}}}
\author{Lea Bold$^1$, Lars Grüne$^2$, Manuel Schaller$^1$, and Karl Worthmann$^1$}
\date{%
	\normalsize
	$^1$Optimization-based Control Group, Technische Universit\"at Ilmenau, Germany\\%
	$^2$Chair of Applied Mathematics, University of Bayreuth, Germany\\[2ex]%
	July 2024
}

\newtheorem{theorem}{Theorem}

\newtheorem{lemma}[theorem]{Lemma}
\newtheorem{proposition}[theorem]{Proposition}
\newtheorem{assumption}[theorem]{Assumption}
\newtheorem{definition}[theorem]{Definition}

\newtheorem{algorithm}[theorem]{Algorithm}

\newcommand{\R}{\mathbb{R}}
\newcommand{\nx}{n_x}
\newcommand{\nc}{n_c}
\newcommand{\N}{M}
\renewcommand{\d}{d}
\newcommand{\calK}{\mathcal K}
\newcommand{\calL}{\mathcal L}
\newcommand{\calU}{\mathcal U}
\newcommand{\bX}{\mathbb X}
\newcommand{\bU}{\mathbb U}

\begin{document}
	
	\maketitle
	
	\noindent \textbf{Abstract}: For nonlinear (control) systems, extended dynamic mode decomposition (EDMD) is a popular method to obtain data-driven surrogate models. 
	Its theoretical foundation is the Koopman framework, in which one propagates observable functions of the state to obtain a linear representation in an infinite-dimensional space. 
	In this work, we prove practical asymptotic stability of a (controlled) equilibrium for EDMD-based model predictive control, in which the optimization step is conducted using the data-based surrogate model. 
	To this end, we derive novel bounds on the estimation error that are proportional to the norm of state and control. 
	This enables us to show that, if the underlying system is cost controllable, this stabilizablility property is preserved. 
	We conduct numerical simulations illustrating the proven practical asymptotic stability.
	
\section{Introduction}

\noindent
Model Predictive Control (MPC; \cite{GrunPann17}) is a well-established feedback control technique. 
In each iteration, an optimal control problem is solved, and a first portion of the optimal control is applied~\cite{CoroGrun20}. This process is then repeated at the successor time instant after measuring (or estimating) the resulting state of the system. 
The popularity of MPC is mainly due to its solid mathematical foundation and the ability to cope with nonlinear constrained multi-input systems. 
In the optimization step, it is, however, necessary to predict the cost functional and/or constraints along the flow of the underlying system, which requires a model, e.g., based on first principles. 

Due to recent progress in data-driven methods, there are several works considering MPC and other model-based controllers using data-driven surrogate models. 
A popular approach is based on extended dynamic mode decomposition (EDMD~\cite{WillKevr15}) as an approximation technique in the Koopman framework. 
The key idea is to lift a nonlinear (control) system to a linear, but infinite-dimensional one and, then, employ EDMD to generate a data-driven finite-dimensional approximation~\cite{OttoRowl21}. 
Convergence of EDMD in the infinite-data limit was shown in~\cite{KordMezi18a}. 
Generally speaking, the Koopman framework can be utilized for data-driven predictions of so-called observables (quantities of interest, e.g., the stage cost in MPC) along the flow of the dynamical (control) system. 
For control systems there are two popular approaches: The first seeks a linear surrogate and is widely called (e)DMDc~\cite{ProcBrun16,KordMezi18b}. 
The second approach yields a bi-linear representation~\cite{WillHema2016} and performs particularly well for systems with direct state-control coupling. 
For this approach also finite-data error bounds for ordinary and stochastic differential equations with i.i.d.\ and ergodic sampling were recently shown in~\cite{SchaWort23,NuskPeit23}. 

In \cite{MamaCast21}, an LQR-based approach to control unconstrained systems by means of a linear surrogate model using Taylor arguments is proposed. 
The performance was further assessed in~\cite{MaHuan19} using a simulation study. 
Recently, robust control of bi-linear Koopman models with guarantees was proposed in~\cite{StraBerb23} or, using Lyapunov-based arguments, in~\cite{SonNara22,NaraSon23}. 
However, without rigorously linking the analysis to verifiable error bounds. 
EDMD-based surrogate models were further applied in the prediction step of MPC~\cite{PeitOtto20,KordMezi18b} and \cite{ZhanPan22} for a robust tube-based approach.
Simulation-based case studies can be found in~\cite{YuShen22} for Koopman-based MPC and in~\cite{KanaYama22} for the bi-linear approach. 
Whereas many of the proposed approaches are shown to perform well in examples, no rigorous guarantees for closed-loop stability of Koopman-based MPC are given.

The main contribution of this work is threefold. 
First, we propose and prove novel error bounds, which are \textit{proportional} to the distance from the desired set point rather than uniform in the state, building upon the error bounds derived in~\cite{NuskPeit23}. 
Second, we show that cost controllability (roughly speaking asymptotic null controllability in terms of the stage costs, see~\cite{CoroGrun20} for details), i.e., a key property to 
establish asymptotic stability in MPC without terminal conditions, is preserved using the EDMD-based surrogat. 
Third, we establish semi-global practical asymptotic stability of the original system if the feedback law is computed using the data-driven surrogate model only. 
To this end, we recall a key result from~\cite{GrunPann17} on practical asymptotic stability for numerical approximations and verify the 
assumptions based on the novel proportional error bounds and the maintained cost controllability. 

The manuscript is organized as follows. In Section~\ref{sec:eDMD}, we recap EDMD within the Koopman framework. 
Then, we introduce MPC, before we derive the novel proportional error bound and provide the problem formulation. 
In Section~\ref{sec:PAS}, we present our main results, i.e., the preservation of cost controllability for the EDMD-based surrogate and practical asymptotic stability of the EDMD-based MPC closed loop. 
Then, we illustrate our findings by means of a simulation study. 
Finally, conclusions are drawn in Section~\ref{sec:conclusions}.

\noindent\textbf{Notation}: 
We use the following comparison functions: $\alpha \in \mathcal{C}(\R_{\geq 0},\R_{\geq 0})$ is said to be of class~$\mathcal{K}$ if it is strictly increasing with $\alpha(0)=0$ and of class~$\mathcal{K}_\infty$ if it, in addition, grows unboundedly. 
A function $\delta \in \mathcal{C}(\R_{\geq 0},\R_{\geq 0})$ is of class~$\mathcal{L}$ if it is strictly decreasing with $\lim_{t\to \infty}\delta(t)=0$. 
Moreover, $\beta \in \mathcal{C}(\R_{\geq 0}^2,\R_{\geq 0})$ is said to be of class~$\mathcal{KL}$ if $\beta(\cdot,t) \in \mathcal{K}$ and $\beta(r,\cdot) \in \mathcal{L}$ hold.
For integers $n \leq m$, we set $[n:m] := [n,m] \cap \mathbb{Z}$. 
The $i$-th standard unit vector in $\mathbb{R}^n$ is denoted by $e_i$, $i\in [1:n]$. For a matrix $A=(a_{ij})\in \R^{n\times m}$, $\|A\|_F^2 = \sum_{i=1}^{n}\sum_{j=1}^m a_{ij}^2$ denotes the squared Frobenius norm. For a set~$X$, we denote the interior by $\operatorname{int}(X)$.

\section{Koopman-based prediction and control}
\label{sec:eDMD}

\noindent 
In this section, we recap the basics of surrogate modeling of nonlinear control systems within the Koopman framework. The underlying idea is to exploit an identity between the nonlinear flow and a linear, but infinite-dimensional operator. 
Then, a compression of this operator onto a finite-dimensional subspace is approximated by extended dynamic mode decomposition (EDMD) using finitely many samples of the system. 

First, we consider the autonomous dynamical system governed by the \emph{nonlinear} ordinary differential equation (ODE)
\begin{align}\label{eq:ode}
\dot{x}(t) = g_0(x(t)),
\end{align}
with locally-Lipschitz continuous map $g_0: \mathbb{R}^{\nx} \rightarrow \mathbb{R}^{\nx}$. 
For initial condition $x(0) = \hat{x} \in \R^{\nx}$, we denote the unique solution of System~\eqref{eq:ode} at time $t \in [0,\infty)$ by $x(t;\hat{x})$. 
We consider the ODE~\eqref{eq:ode} on a compact and non-empty set $\bX \subset \mathbb{R}^{\nx}$. Then, to avoid technical difficulties in this introductory section, forward invariance of the set~$\bX$ w.r.t.\ the dynamics~\eqref{eq:ode} is assumed, i.e., $x(t;\hat{x}) \in \mathbb{X}$, $t \geq 0$, holds for all $\hat{x} \in \mathbb{X}$. This may be ensured, e.g., by some inward-pointing condition and guarantees existence of the solution on~$[0,\infty)$. 
Then, the Koopman semigroup $(\mathcal{K}^t)_{t \geq 0}$ of bounded linear operators is defined by the identity 
\begin{equation}\label{eq:Koopman}
(\mathcal{K}^t \varphi)(\hat{x}) = \varphi(x(t;\hat{x})) \quad\forall\,t \geq 0, \hat{x} \in \bX, \varphi \in L^2(\bX,\mathbb{R}),
\end{equation}
see, e.g., \cite[Prop.~2.4]{MaurSusu20} or \cite[Chapter 7]{LasoMack13}. Here, the real-valued functions~$\varphi$ are called \emph{observables}. 
The identity~\eqref{eq:Koopman} states that, instead of evaluating the observable~$\varphi$ at the solution of the \emph{nonlinear} 
system~\eqref{eq:ode} emanating from initial state~$\hat{x}$ at time~$t$, one may also apply the linear, \emph{infinite-dimensional} Koopman operator~$\mathcal{K}^t$ to the observable~$\varphi$ and, then, evaluate $\mathcal{K}^t \varphi$ at~$\hat{x}$.

Since the flow of System~\eqref{eq:ode} is continuous, $(\calK^t)_{t \geq 0}$ is a strongly-continuous semigroup of bounded linear operators. 
Correspondingly, we can define the, in general, unbounded infinitesimal generator~$\mathcal{L}$ of this semigroup by
\begin{align}\label{eq:generator}
\calL \varphi := \lim_{t \searrow 0}\frac{\calK^t \varphi - \varphi}{t} \qquad \forall\, \varphi \in D(\calL),
\end{align}
where the domain $D(\calL)$ consists of all $L^2$-functions, for which the above limit exists. 
Using this generator, we may formulate the equivalent evolution equation for $\Phi(t)= \calK^t\varphi = \varphi(x(t;\cdot))$ 
\begin{align}\label{eq:Cauchy_problem}
\dot \Phi(t) = \calL \Phi(t), \qquad \Phi(0)=\varphi.
\end{align}

\noindent Next, we recap the extension of the Koopman approach to control-affine systems, i.e., systems governed by the dynamics
\begin{equation}\label{eq:dynamics_control_affine}
\dot{x}(t) = g_0(x(t)) + \sum_{i=1}^{\nc} g_i(x(t)) u_i(t),
\end{equation}
where the control function $u \in L^\infty_{\operatorname{loc}}([0,\infty),\R^{\nc})$ serves as an input and the input maps $g_i:\R^{\nx} \to \R^{\nx}$, $i \in [0:\nc]$, are locally Lipschitz continuous.
A popular approach to obtain a data-based surrogate model is DMDc~\cite{ProcBrun16} or c~\cite{KordMezi18b}, where one seeks a linear control system. In this paper, we pursue an alternative \emph{bi-linear} approach, which exploits the control-affine structure of 
system~\eqref{eq:dynamics_control_affine} and was --~to the best of our knowledge~-- proposed by~\cite{WillHema2016,Sura16}. 
This approach shows a superior performance for systems with state-control coupling~\cite{BrudFu21,FolkBurd21}. 
For the flow of the control system \eqref{eq:dynamics_control_affine} with constant control input $u$, the Koopman operator $\calK^t_u$ is defined analogously to~\eqref{eq:Koopman}.
A straightforward computation shows that its generator preserves control affinity, i.e.,
\begin{align}\label{eq:generator_control}
\calL^{u} = \calL^0 + \sum\nolimits_{i=1}^{\nc} u_i(\calL^{e_i}-\calL^0)
\end{align}
holds for $u \in \mathbb{R}^{\nc}$,
where $\calL^0$ and $\calL^{e_i}$, $i\in [1:\nc]$, are the generators of the Koopman semigroups corresponding to the constant controls $u \equiv 0$ and $u \equiv e_i$, $i \in [1:\nc]$, respectively.
For general control functions~$u \in L^\infty_{\operatorname{loc}}([0,\infty),\R^{\nc})$, one can now state the respective abstract Cauchy problem analogously to~\eqref{eq:Cauchy_problem} replacing the generator~$\mathcal{L}$ by its time-varying counterpart~$\mathcal{L}^{u(t)}$ defined by~\eqref{eq:generator_control}, see~\cite{NuskPeit23} for details.

The success of the Koopman approach in recent years is due to its linear nature such that the compression of the Koopman operator or its generator~\eqref{eq:generator_control} to a finite-dimensional subspace -- called dictionary -- leads to matrix representations. 
Being finite-dimensional objects, these matrices can then be approximated by a finite amount of data.
Let the dictionary $\mathbb{V} := \operatorname{span}(\{ \psi_k: k \in [1:\N]\} )$ be the $\N$-dimensional subspace spanned by the chosen observables~$\psi_k$. 
We denote the $L^2$-orthogonal projection onto $\mathbb{V}$ by $P_\mathbb{V}$. Further, using $\d$ i.i.d.\ data points $x_1,\ldots,x_{\d}\in \bX$, the $(\N \times \d)$-matrices
\begin{align*}
X := \left( \left.   \left(\begin{smallmatrix}
\psi_1(x_1)\\
:\\
\psi_{\N}(x_1)
\end{smallmatrix}\right)\right| \ldots \left| \left(\begin{smallmatrix}
\psi_1(x_{\d})\\
:\\
\psi_{\N}(x_{\d})
\end{smallmatrix}\right)\right. \right)\quad\text{ and }\quad
Y := \left( \left. \left(\begin{smallmatrix}
(\mathcal{L}^0\psi_1)(x_1)\\
:\\
(\mathcal{L}^0\psi_{\N})(x_1)
\end{smallmatrix}\right)\right| \ldots \left| \left(\begin{smallmatrix}
(\mathcal{L}^0\psi_1)(x_{\d})\\
:\\
(\mathcal{L}^0\psi_{\N})(x_{\d})
\end{smallmatrix}\right)\right. \right)
\end{align*}
are defined, where $(\mathcal{L}^0 \psi_k)(x_j) = \nabla \psi_k(x_j)^\top g_0(x_j)$ holds for $k \in [1:\N]$ and $j \in [1:\d]$. 
Then, the empirical estimator of the compressed Koopman generator $P_\mathbb{V}\calL^0\vert_\mathbb{V}$ is given by
\begin{align*}
{\calL}^0_\d := \operatorname{arg}\min_{{L} \in \R^{\N\times \N}} \|{L}X-Y\|_F^2.
\end{align*}
We have to repeat this step for $\mathcal{L}^{e_i}$, $i \in [1:\nc]$, based on the identity $$(\mathcal{L}^{e_i} \psi_k)(x_j) = \nabla \psi_k(x_j)^\top \left(g_0(x_j) + g_i(x_j)\right)$$ to construct the data-driven approximation of~$\mathcal{L}^{u}$ according to~\eqref{eq:generator_control}. 
Consequently, for $\varphi \in \mathbb{V}$ and control function $u \in L^\infty_{\operatorname{loc}}([0,t],\R^{\nc})$, a data-driven predictor is given as the solution of the linear time-varying Cauchy problem~\eqref{eq:Cauchy_problem}, where the unbounded operator~$\mathcal{L}$ is replaced by $\calL_d^{u(t)}$. 
The convergence of this estimator was shown in~\cite{KordMezi18a} if both the dictionary size and the number of data points goes to infinity. 
Finite-data bounds typically split the error into two sources: A projection error stemming from the finite dictionary and an estimation error resulting from a finite amount of data. 
A bound on the estimation error for control systems was derived in~\cite{NuskPeit23}, where, in addition to i.i.d.\ sampling of ODEs, also SDEs and ergodic sampling, i.e.\ sampling along one sufficiently-long trajectory, were considered. 
A full approximation error bound for control systems was provided in~\cite{SchaWort23} using a dictionary of finite elements. 
We provide an error bound tailored to the sampled-data setting used in this work in Subsection~\ref{subsec:sampled}.

\section{Proportional error bound for EDMD-based MPC and problem formulation}
\label{sec:problem}

\noindent We consider the discrete-time control system given by
\begin{equation}\label{eq:dynamics_DT}
x^+ = f(x,u)    
\end{equation}
with nonlinear map $f: \mathbb{R}^{\nx} \times \mathbb{R}^{\nc} \rightarrow \mathbb{R}^{\nx}$.
Then, for initial state $\hat{x} \in \mathbb{R}^{\nx}$ and sequence of control values $(u(k))_{k \in \mathbb{N}_0}$, $x_u(n;\hat{x})$ denotes the solution at time~$n \in \mathbb{N}_0$, which is recursively defined by~\eqref{eq:dynamics_DT} and $x_u(0;\hat{x}) = \hat{x}$. 
In the following, $f(0,0) = 0$ is assumed, i.e., the origin is a controlled equilibrium for~$u=0$. After reviewing the basics of model predictive control, we derive a sampled-data representation of the continuous-time dynamics~\eqref{eq:dynamics_control_affine} and the corresponding abstract Cauchy problem, i.e., \eqref{eq:Cauchy_problem} with $\calL^{u(t)}$ including its -based surrogate in Subsection~\ref{subsec:sampled}. 
Then, we provide the problem formulation in Subsection~\ref{subsec:problem}.

\noindent We impose state and control constraints using the compact sets $\bX \subset \mathbb{R}^{\nx}$ and $\bU \subset \mathbb{R}^{\nc}$ with~$(0,0) \in \operatorname{int}(\bX \times \bU)$, respectively. Next, we define admissibility of a sequence of control values. \begin{definition}\label{def:admissibility}
	A sequence of control values $(u(k))_{k=0}^{N-1} \subset~\bU$ of length~$N$ is said to be admissible for state $\hat{x} \in \bX$, if $x_u(k;\hat{x}) \in \bX$ holds for all $k \in [1:N]$. 
	For $\hat{x} \in \bX$, the set of admissible control sequences is denoted by~$\calU_N(\hat{x})$. If, for $u = (u(k))_{k \in \mathbb{N}_0}$, $(u(k))_{k=0}^{N-1} \in \calU_N(\hat{x})$ holds for the restriction of~$u$ for all $N \in \mathbb{N}_0$, we write $u \in \calU_\infty(\hat{x})$.
\end{definition}

We introduce the quadratic stage cost $\ell: \bX \times \bU \rightarrow \mathbb{R}_{\geq 0}$, 
\begin{equation}\label{eq:stage_cost}
\ell(x,u) := \| x \|_Q^2 + \| u \|_R^2 := x^\top Q x + u^\top R u,
\end{equation}
for symmetric and positive definite matrices $Q \in \mathbb{R}^{\nx \times \nx}$ and $R \in \mathbb{R}^{\nc \times \nc}$. Next, based on Definition~\ref{def:admissibility}, we introduce the MPC Algorithm, where we tacitly assume existence of an optimal sequence of control values in Step~(2) along the MPC closed-loop dynamics and full-state measurement.
\begin{algorithm}[Model Predictive Control with horizon~$N$]\label{alg:MPC}
	At each time~$n \in \mathbb{N}_0$: 
	\begin{enumerate}
		\item [(1)] Measure the state $x(n) \in \mathbb{X}$ and set $\hat{x} := x(n)$. 
		\item [(2)] Solve the optimization problem
		\begin{equation}\nonumber
		u^\star \!\in\! \operatorname{argmin}_{u \in \calU_N(\hat{x})}\ J_N(\hat{x},u) :=\! \sum_{k=0}^{N - 1} \ell(x_u(k;\hat{x}),u(k)) 
		\end{equation}
		subject to $x_u(0;\hat{x}) = \hat{x}$ and the dynamics $x_u(k+1;\hat{x}) = f(x_u(k;\hat{x}),u(k))$, $k \in [0:N-2]$.
		\item [(3)] Apply the feedback value $\mu_N(x(n)) := u^\star(0) \in \bU$.\\ 
	\end{enumerate}
\end{algorithm}

\noindent Overall, Algorithm~\ref{alg:MPC} yields the MPC closed-loop dynamics
\begin{equation}\label{eq:dynamics_closed_loop}
x^+_{\mu_N} = f(x_{\mu_N}, \mu_N(x_{\mu_N})),
\end{equation}
where the feedback law~$\mu_N$ is well defined at~$\hat{x}$ if $\calU_N(\hat{x}) \neq \emptyset$ holds. 
We emphasize that this condition holds if, e.g., $\bX$ is controlled forward invariant and refer to~\cite{BoccGrun14} and~\cite{EsteWort20} for sufficient condition to ensure recursive feasibility without requiring controlled forward invariance of $\bX$ (and without terminal conditions) for discrete and continuous-time systems, respectively. 
The closed-loop solution resulting from the dynamics~\eqref{eq:dynamics_closed_loop} is denoted by $x_{\mu_N}(n;\hat{x})$, where $x_{\mu_N}(0;\hat{x}) = \hat{x}$ holds. 
Moreover, we define the (optimal) value function~$V_N: \bX \rightarrow \mathbb{R}_{\geq 0} \cup \{ \infty \}$ as $V_N(x) := \inf_{u \in \calU_N(x)} J_N(x,u)$.

\subsection{Proportional error bound for sampled-data systems}
\label{subsec:sampled}

\noindent 
We consider the nonlinear continuous-time control system given by~\eqref{eq:dynamics_control_affine}. 
Equidistantly discretizing the time axis~$[0,\infty)$, i.e., using the partition $\bigcup_{k=0}^\infty [k \Delta t, (k+1) \Delta t)$ with sampling period~$\Delta t > 0$, and using a (piecewise) constant control function on each sampling interval, i.e., $u(t) \equiv \hat{u} \in \bU \subset \mathbb{R}^{\nc}$ on~$[k \Delta t, (k+1) \Delta)$, we generate the discrete-time system
\begin{equation}\label{eq:dynamics_sampled_data}
x^+ \hspace*{-0.5mm}= f(\hat{x},\hat{u}) := \!\int_0^{\Delta t}\hspace*{-0.2cm} g_0(x(t;\hat{x},u)) + \sum_{i=1}^{\nc} g_i(x(t;\hat{x},u)) u_i(t)\,\mathrm{d}t.
\end{equation}
We emphasize that the drift~$g_0$ does not exhibit an offset independently of the state variable~$x$ in view of our assumption $f(0,0) = 0 = g_0(0)$.
We define the vector-valued observable
\begin{align}\label{eq:dictionary}
\begin{split}
\Psi(x) &= \begin{pmatrix}
\psi_1(x),\ldots,\psi_\N(x)
\end{pmatrix} \\&= \begin{pmatrix}
1, x_1,\ldots,x_{\nx},\psi_{\nx+2}(x),\ldots,\psi_\N(x)
\end{pmatrix},
\end{split}
\end{align}
where $\psi_1(x)\equiv 1$, $\psi_{k+1}(x) = x_k$, $k \in [1:\nx]$, and $\psi_k \in \mathcal{C}^1(\mathbb{R}^{\nx},\mathbb{R})$, $k \in [\nx + 2:\N]$, are locally-Lipschitz continuous functions satisfying $\psi_k(0) = 0$ and $(D \psi_k)(0) = 0$. 
Hence, $\Psi:\bX \to \mathbb{R}^\N$ is Lipschitz continuous with 
constant~$L_\Psi$ such that $\| \Psi(x) - \Psi(0) \| \leq L_\Psi \| x \|$ holds. 
A straightforward calculation then shows $(P_{\mathbb{V}} \mathcal{L}^0|_{\mathbb{V}})_{k,1} \equiv 0$, $k \in [1:\N]$, which we impose for the data-driven approximation to ensure consistency, i.e., that $f(0,0) = g_0(0) = 0$ is preserved. For $g_i$, $i \in [1:n_c]$, the first (constant) observable enables us to approximate components of the control maps, which do not depend on the state~$x$, separately.

In this note, we make use of the following Assumption~\ref{ass:invariance}, which ensures that no projection error occurs. This assumption is common in systems and control when the Koopman framework is used, see, e.g., \cite{ProcBrun18,KordMezi18b}. The construction of suitable dictionaries ensuring this assumption is discussed in \cite{BrunBrun16,KordMezi20}. 
A condition ensuring this invariance is provided, e.g., in~\cite[Theorem~1]{GoswPale21}, where even a method for the construction of a suitable dictionary is discussed. 
\begin{assumption}[Invariance of~$\mathbb{V}$] 
	\label{ass:invariance}
	For any $\varphi \in \mathbb{V}$, the relation $\varphi(x(\Delta t;\cdot,u)) \in \mathbb{V}$ holds for all $u(t) \equiv \hat{u} \in \bU \subset \R^{\nc}$.
\end{assumption}

We note that if this invariance assumption does not hold, and in order mitigate the projection error, subspace identification methods may be employed to (approximately) ensure invariance of the dictionary, i.e., the space spanned by the choosen observables, see, e.g., \cite{HaseCort21,KrolTell22}.

Next, we deduce an error bound adapted to our sampled-data setting. Assumption~\ref{ass:invariance} implies that the compression of the generator coincides with its restriction onto~$\mathbb{V}$, i.e., $P_\mathbb{V} \mathcal{L}^u\vert_\mathbb{V} = \mathcal{L}^u\vert_\mathbb{V}$. Thus, for $u\in \bU$, the Koopman operator is the matrix exponential of the generator, i.e., $\calK^{\Delta t}_u = e^{\Delta t \mathcal{L}^u}$ holds.
\begin{proposition}\label{prop:generatorbound}
	Suppose that Assumption~\ref{ass:invariance} holds. 
	For every error bound $\varepsilon>0$ and 
	probabilistic tolerance $\delta \in (0,1)$, there is an amount of data $\d_0\in \mathbb{N}$ such that with probability $1-\delta$, the error bound
	\begin{align}\label{eq:operatorbound}
	\big\|e^{\Delta t \mathcal{L}^u \vert_{\mathbb{V}}} - e^{\Delta t \calL_\d^u}\big\| \leq \varepsilon
	\end{align}
	holds for all $\d \geq \d_0$ and all $u\in \bU$ for the Koopman operator $\calK^{\Delta t}_u = e^{\Delta t \mathcal{L}^u}$.
\end{proposition}
\begin{proof}
	For $g(t) = e^{t \calL^u\vert_{\mathbb{V}}}-e^{t \calL_\d^u}$, we have
	\begin{align*}
	g^\prime(t)  = \calL^u\vert_{\mathbb{V}} e^{t \calL^u\vert_{\mathbb{V}}} \mp \calL^u\vert_{\mathbb{V}}e^{t \calL_\d^u} - \calL_\d^ue^{t \calL_\d^u} 
	= \calL^u\vert_{\mathbb{V}}g(t) + (\calL^u\vert_{\mathbb{V}}-\calL_\d^u) \Big( e^{t \calL_\d^u} \mp e^{t \calL^u\vert_{\mathbb{V}}} \Big).
	\end{align*}
	Since $g(0)=0$, we have $g(t) = \int_0^{\Delta t} g^\prime(s)\,\mathrm{d}s$. Then, plugging in the 
	expression for $g^\prime(s)$, the triangle inequality yields
	\begin{align*}
	& \|g(t)\| \leq \int_0^{t} \beta \|g(s)\|\,\mathrm{d}s + \alpha(t) 
	\end{align*}
	with the constant $\beta = \|\calL^u_*\vert_{\mathbb{V}}\| + \|(\calL^u\vert_{\mathbb{V}}-\calL_\d^u)\|$ and
	\begin{align*}\nonumber
	\alpha(t) = \|(\calL^u\vert_{\mathbb{V}} - \calL_\d^u) \| \int_0^{t} \|e^{s \calL^u\vert_{\mathbb{V}}} \|\,\mathrm{d}s 
	 \leq \frac {\Delta t \cdot \|(\calL^u\vert_{\mathbb{V}} - \calL_\d^u) \|}{\| \calL^u_*\vert_{\mathbb{V}} \|} \Big( e^{\Delta t \| \calL^u_*\vert_{\mathbb{V}} \|} - 1 \Big) =: c_{\Delta t}
	\end{align*}
	for all $t \in (0,\Delta t]$, where $\calL^u_*\vert_{\mathbb{V}}$ maximizes $\| \calL^u\vert_{\mathbb{V}} \|$ w.r.t.\ the compact set~$\mathbb{U}$.
	Then, Gronwall's inequality with $\alpha(t)$ replaced by $c_{\Delta t}$ yields
	\begin{align*}
	\| g(\Delta t) \| & \leq c_{\Delta t} \Big( 1 + \int_0^{\Delta t} \beta e^{(\Delta t-t)\beta}\,\mathrm{d}t \Big) = c_{\Delta t} e^{\Delta t \beta}. 
	\end{align*}
	Invoking \cite[Theorem 3]{SchaWort23} yields, for any $\tilde \varepsilon> 0$, a sufficient amount of data $d_0\in \mathbb{N}$ such that $\|\mathcal{L}^u\vert_\mathbb{V} -\mathcal{L}^u_\d \| \leq \tilde \varepsilon$ holds for all $u \in \mathbb{U}$ and $d \geq d_0$. 
	Hence, setting $\tilde{\varepsilon}$ such that the inequality
	\begin{equation}
	\frac {\Delta t \cdot \tilde{\varepsilon}}{\| \calL^u_*\vert_{\mathbb{V}} \|} \Big( e^{\Delta t \| \calL^u_*\vert_{\mathbb{V}} \|} - 1 \Big) e^{\Delta t \left( \|\calL^u_*\vert_{\mathbb{V}}\| + \tilde{\varepsilon} \right) } \leq \varepsilon
	\end{equation}
	holds and using the definitions of~$\beta$ and $c_{\Delta t}$ ensures Inequality~\eqref{eq:operatorbound}. 
	Since the left hand side is monotonically increasing in $\tilde{\varepsilon}$ and zero for $\tilde{\varepsilon} = 0$, this is always possible, which completes the proof.
\end{proof}

We briefly quantify the sufficient amount of data~$d_0$ 
in view of the dictionary size~$\N$ and the parameters~$\varepsilon$ and~$\delta$. 
First, by a standard Chebychev inequality, one obtains the dependency $d_0 \sim \nicefrac{\N^2}{\varepsilon^2\delta}$, cf.\ \cite{SchaWort23,NuskPeit23}. 
This can be improved in reproducing kernel Hilbert spaces, where the dictionary is given by feature maps given by the kernel evaluated at the samples. Here a scaling depending logarithmically on $\delta$ was shown in \cite[Proposition 3.4]{PhilScha23} using Hoeffding's inequality, see also~\cite{PhilScha2024}. In the latter reference, invariance conditions were discussed, which may allow to relax Assumption~\ref{ass:invariance}. 
Otherwise, only bounds on the projection error w.r.t.\ the $L_2$-norm are available~\cite{SchaWort23}, which does not yield pointwise bounds.

For the discrete-time dynamics~\eqref{eq:dynamics_sampled_data}, we get the identity
\begin{align}\label{eq:dynamics_exact}
f(\hat{x},\hat{u}) = P_{x} e^{\Delta t \mathcal{L}^{\hat{u}}|_{\mathbb{V}}} \Psi(\hat{x})
\end{align}
resulting from sampling with zero-order hold in view of Assumption~\ref{ass:invariance}, where $P_x:\mathbb{R}^\N \to \mathbb{R}^{\nx}$ is the projection onto the first $\nx$~components.
Further, based on the bi-linear -based surrogate model of Subsection~\ref{sec:eDMD} for $\d$~data points, we define the data-driven surrogate model
\begin{align}\label{eq:dynamics_approx}
f^\varepsilon(\hat{x},\hat{u}) = P_xe^{\Delta t \mathcal{L}^{\hat{u}}_\d} \Psi(\hat{x}).
\end{align}
Next, we derive a novel error bound that is \textit{proportional} to the norm of the state and the control and, thus, ensures that the error becomes small close to the origin.
\begin{proposition}\label{prop:errbound}
	Let~$L_\Psi$ be the Lipschitz constant of $\Psi$ on the set~$\bX$. 
	Then, for every 
	error bound $\varepsilon \in (0,\varepsilon_0]$, the inequality
	\begin{align}\label{eq:dynamics_bound}
	\|f(x,u)-f^\varepsilon(x,u)\| \leq \varepsilon \left( L_\Psi \| x \| + \Delta t \cdot \tilde{c} \| u \| \right)
	\end{align}
	holds for all $x \in \bX$ and $u \in \bU$ with some constant~$\tilde{c}$ if \eqref{eq:operatorbound} holds provided $\{ f(x,u), f^\varepsilon(x,u) \} \subset \bX$.
\end{proposition}
\begin{proof}
	By local Lipschitz continuity of $\Psi$, $0\in \operatorname{int}(\bX)$ and $\| P_x \| \leq 1$ we compute
	\begin{align*}
	 \|f(x,u) \hspace*{-0.5mm}-\hspace*{-0.5mm} f^\varepsilon(x,u)\| &= \big\|P_x [e^{\Delta t \mathcal{L}^u\vert_\mathbb{V}} - e^{\Delta t \mathcal{L}^u_\d}] [ \Psi(x) \pm \Psi(0)] 
	\big\| \\ 
	&\leq \underbrace{\varepsilon \|\Psi(x) - \Psi(0) \|}_{\leq L_\Psi \varepsilon \| x \|} + \big\| \underbrace{ ( e^{\Delta t \mathcal{L}^u \vert_\mathbb{V}} - e^{\Delta t \calL_\d^u} ) \Psi(0) }_{ =:h(\Delta t) } \big\|. 
	\end{align*}
	Then, Taylor series 
	expansion of $h(\Delta t) = h(0) + \Delta t \cdot h^\prime(\xi)$, $\xi \in [0,\Delta t]$, with $h(0) = 0$ leads to the representation
	\begin{align*}
	\text{\scriptsize $\frac {h(\Delta t)}{\Delta t}$} = & ( e^{\xi \mathcal{L}^u \vert_\mathbb{V}} \mathcal{L}^u \vert_\mathbb{V} \pm e^{\xi \calL_\d^u} \mathcal{L}^u \vert_\mathbb{V} - e^{\xi \calL_\d^u} \calL_\d^u ) \Psi(0) \\
	= & ( e^{\xi \mathcal{L}^u \vert_\mathbb{V}} - e^{\xi \calL_\d^u} ) \mathcal{L}^u \vert_\mathbb{V} \Psi(0) + e^{\xi \calL_\d^u} \left( \mathcal{L}^u \vert_\mathbb{V} - \calL_\d^u \right) \Psi(0).
	\end{align*}%
	For a sufficient amount of data~$d_0 \in \mathbb{N}$, we have $\max_{i \in [1:n]} \| \mathcal{L}^{e_i}|_{\mathbb{V}} - \mathcal{L}^{e_i}_d \| \leq \bar{\varepsilon}$. Then, the second summand can be estimated by
	\begin{align}\nonumber 
	\left[ \| e^{\xi \mathcal{L}^u \vert_\mathbb{V}} \hspace*{-0.5mm}-\hspace*{-0.5mm} e^{\xi \calL_\d^u} \| + \| e^{\xi \mathcal{L}^u \vert_\mathbb{V}} \| \right] \hspace*{-.5mm} \left\| \left( \mathcal{L}^u \vert_\mathbb{V} - \calL_\d^u \right) \Psi(0) \right\| \leq c_0 \bar{\varepsilon} \| u \|
	\end{align}
	with $c_0 := e^{\Delta t \| \mathcal{L}^u_*|_{\mathbb{V}} \|} + \varepsilon$ with $\varepsilon$ from Proposition~\ref{prop:generatorbound}, 
	where $\calL^u_*\vert_{\mathbb{V}}$ maximizes 
	$\| \mathcal{L}^u\vert_{\mathbb{V}} \|$ w.r.t.\ the compact set~$\mathbb{U}$ and we have used that the contributions of~$\mathcal{L}^0$ and~$\mathcal{L}^0_d$ cancel out thanks to $\Psi(0)$ and the control value acts as a factor. 
	The same argument yields 
	$\| \mathcal{L}^u \vert_\mathbb{V} \Psi(0) \| \leq \| \mathcal{L}^u_*|_{\mathbb{V}} \| \| u \|$. 
	Combining the derived estimates yields the assertion, i.e., Inequality~\eqref{eq:dynamics_bound} with $\tilde{c} := e^{\Delta t \| \mathcal{L}^u_*|_{\mathbb{V}} \|} + \varepsilon_0 + \| \mathcal{L}^u_*|_{\mathbb{V}} \|$.
\end{proof}

In~\cite{StraBerb23}, a bound of the form~\eqref{eq:dynamics_bound} was assumed in the lifted space, i.e., without the projector~$P_x$. Therein, the bound was used to construct a feedback controller achieving robust local stability using a finite gain argument. However, the bound was not established, but rather assumed --~in addition to the invariance in Assumption~\ref{ass:invariance}.

\subsection{Problem statement}
\label{subsec:problem}

\noindent 
We will leverage the error bound of Proposition~\ref{prop:errbound} to provide a stability result when using the surrogate dynamics~$f^\varepsilon$ in Step (2) of the MPC Algorithm~\ref{alg:MPC} to stabilize the original system. 
The main result shows that, if the nominal MPC controller is asymptotically stabilizing, the data-based controller with $f^\varepsilon$ ensures convergence to a neighborhood of the origin, whose size depends on~$\varepsilon$, i.e., practical asymptotic stability. 
\begin{definition}[Practical asymptotic stability]\label{def:pracstab}
	For $\varepsilon>0$, let $\mu_N^\varepsilon$ be the feedback law defined in Algorithm~\ref{alg:MPC} with $f=f^\varepsilon$, where admissibility of control sequences at $\hat{x}$, i.e., $u \in \calU_N^\varepsilon(\hat{x})$, is defined w.r.t.\ the tightened set $\bX \ominus \mathcal{B}_\varepsilon(0)$. 
	Let $A \subset \mathbb{X} \ominus \mathcal{B}_\varepsilon(0)$ be given such that $\mathcal{U}_N^\varepsilon(\hat{x}) \neq \emptyset$ for all $\hat{x} \in A$.
	Then, the origin is said to be semi-globally practically asymptotically stable (PAS) on~$A$ 
	if there exists $\beta \in \mathcal{K}\mathcal{L}$ such that for each $r>0$ and $R > r$ there is $ \varepsilon_0 > 0$ such that for each $\hat{x} \in A$ with $\| \hat{x} \|\leq R$ and all $\varepsilon \in (0,\varepsilon_0]$ such that \eqref{eq:dynamics_bound} holds, the solution $x_{\mu_N^\varepsilon}(\cdot,\hat{x})$ of 
	\begin{align}\label{eq:ex_cl}
	x_{\mu_N^\varepsilon} (n+1) = f(x_{\mu_N^\varepsilon}(n),\mu_N^\varepsilon(x_{\mu_N^\varepsilon}(n)))
	\end{align}
	with $x_{\mu_N^\varepsilon}(0) = \hat{x}$ satisfies $x_{\mu_N^\varepsilon}(n;\hat{x})\in A$ and 
	\begin{align*}
	\| x_{\mu_N^\varepsilon}(n;\hat{x}) \| \leq \max \{\beta(\| \hat{x} \|,n),r\} \qquad\forall\,n \in \mathbb{N}_0.
	\end{align*}
\end{definition}
The incorporation of the Pontryagin difference $\mathbb{X} \ominus \mathcal{B}_\varepsilon(0)$ in the admissibility of control sequences for the surrogate model ensures that the original system evolves in the compact set~$\bX$, i.e., that every optimal control function is, in particular, admissible for the original system in view of the error bound of Proposition~\ref{prop:generatorbound}.
In the following section, we will show that the error bound shown in Proposition~\ref{prop:errbound} and cost-controllability of the original dynamics imply practical asymptotic stability of the closed-loop using EDMD-based MPC.

\section{Practical asymptotic stability of surrogate-based MPC}
\label{sec:PAS}

\noindent
In this section, we prove our main result, i.e., practical asymptotic stability of the data-based MPC Algorithm~\ref{alg:MPC} using the  surrogate $f^\varepsilon$ as defined in \eqref{eq:dynamics_approx} to stabilize the original system with $f$ given by~\eqref{eq:dynamics_sampled_data} or, equivalently, \eqref{eq:dynamics_exact}.

We follow the line of reasoning outlined in~\cite[Section~11.5]{GrunPann17}. To this end, we recall \cite[Theorem 11.10]{GrunPann17} regarding \emph{stability for perturbed solutions} in Proposition \ref{thm:PAS}, which is a key tool for our analysis. We define
$$
V^\varepsilon_N(\hat{x}) := \inf_{u \in \calU_N^\varepsilon(\hat{x})}\sum_{k=0}^{N - 1} \ell(x^\varepsilon_u(k;\hat{x}),u(k))
$$
where $x_u^\varepsilon(0;\hat{x}) = \hat{x}$ and $x_u^\varepsilon(k+1;\hat{x}) = f^\varepsilon(x_u^\varepsilon(k;\hat{x}),u(k))$ for $k\in [0:N-2]$. 
\begin{proposition}\label{thm:PAS}
	Consider the MPC-feedback law~$\mu_N^\varepsilon$ of Algorithm~\ref{alg:MPC} with $f = f^\varepsilon$, where $f^{\varepsilon}$ satisfies Condition~\eqref{eq:dynamics_bound} and let $S \subset \bX$ be a forward-invariant set w.r.t.\ $f^\varepsilon(\cdot,\mu^\varepsilon_N(\cdot))$.    Further, let the following assumptions hold:
	
	(i) There is $\varepsilon_0>0$ and $\alpha \in (0,1]$ such that for all $\varepsilon\in (0,\varepsilon_0]$ the relaxed dynamic programming inequality
	\begin{align*}
	V_N^\varepsilon(x) \geq \alpha \ell(x,\mu^\varepsilon_N(x)) + V_N^\varepsilon(f^\varepsilon(x,\mu^\varepsilon_N(x)))
	\end{align*}
	holds on $S$. 
	In addition, there exist $\alpha_1,\alpha_2,\alpha_3 \in \mathcal{K}_\infty$ such that
	\begin{align*}
	\alpha_1(\|x\|) \leq V_N^\varepsilon(x) \leq \alpha_2(\|x\|) \ \ \mathrm{and}\ \  \ell(x,u) \geq \alpha_3(\|x\|)
	\end{align*}
	hold for all $x \in S$, $\varepsilon \in (0,\varepsilon_0]$, and $u\in \mathbb{U}$.
	
	(ii) $V^\varepsilon_N$ is uniformly continuous and $f^\varepsilon$ is uniformly continuous in~$u$ on closed balls $\overline{B}_\rho(0)$, i.e., there is $\varepsilon_0$ such that, for each $\rho >0$, there exists $\omega_V, \omega_f \in \mathcal{K}$: 
	\begin{align*}
	|V^\varepsilon_N(x)-V^\varepsilon_N(y)| & \leq \omega_V(\|x-y\|), \\ 
	\|f^\varepsilon(x,u)-f^\varepsilon(y,u)\| & \leq \omega_f(\|x-y\|) \qquad\forall\, u \in \mathbb{U}
	\end{align*}
	for all $x,y \in \overline{B}_\rho(0) \cap S$ and $\varepsilon\in (0,\varepsilon_0]$.
	Then the exact closed-loop system with perturbed feedback~$\mu_N^\varepsilon$ defined in \eqref{eq:ex_cl} is semiglobally practically asymptotically stable on $A = S$ in the sense of Definition~\ref{def:pracstab}.
\end{proposition}

We first verify the condition of Proposition~\ref{thm:PAS} considering uniform continuity of the surrogate model.
\begin{lemma}
	\label{lem:as12}
	Let $\varepsilon_0 > 0$ be given. Then, $f^\varepsilon$ is \textit{uniform continuous} in~$u$ with $\omega_f(r) = c L_\Psi r$, $c = c(\varepsilon_0)$, i.e., 
	\begin{align}\label{eq:feps_lipschitz}
	\|f^\varepsilon(x,u)-f^\varepsilon(y,u)\| \leq c L_\Psi \|x-y\|
	\end{align}
	holds for all $x,y \in \bX$, $u\in \mathbb{U}$, and $\varepsilon \in (0,\varepsilon_0]$ provided that the error bound~\eqref{eq:operatorbound} is satisfied.
\end{lemma}
\begin{proof}
	The error bound~\eqref{eq:operatorbound} and $\|P_x\| \leq 1$ imply
	\begin{align*}
	\|f^\varepsilon(x,u) - f^\varepsilon(y,u)\| & = \| P_xe^{\Delta t \mathcal{L}^u_\d}(\Psi(x)-\Psi(y)) \| \\
	& \leq \|e^{\Delta t \mathcal{L}^u_\d} \mp e^{\Delta t \mathcal{L}^u|_{\mathbb{V}}} \| \hspace*{-0.05cm}\cdot\hspace*{-0.05cm} \|\Psi(x)-\Psi(y)\| \\
	& \leq ( \varepsilon_0 + \| e^{\Delta t \mathcal{L}^u_*|_{\mathbb{V}}} \|) L_\Psi \| x - y \|,
	\end{align*}
	where $\calL^u_*\vert_{\mathbb{V}}$ maximizes $\| \mathcal{L}^u\vert_{\mathbb{V}} \|$ w.r.t.\ the compact set~$\mathbb{U}$. This completes the proof with $c := \varepsilon_0 + \| e^{\Delta t \mathcal{L}^u_*|_{\mathbb{V}}} \|$.
\end{proof}

Using the novel proportional error bound of Proposition~\ref{prop:errbound} we rigorously show that cost controllability as defined in~\cite{CoroGrun20} and~\cite{Wort11} for continuous- and discrete-time systems, respectively, is inherited by the EDMD-based surrogate model. Cost controllability links stabilizability with the stage cost employed in MPC, see, e.g., \cite{GrunPann10,Wort11}.
The only additional requirement is that optimal control sequences have to be admissible also for the surrogate model. While this may be a severe restriction close to the boundary of the set~$\bX \ominus \mathcal{B}_\varepsilon(0)$, it is typically satisfied on a suitably chosen sub-level set of the optimal value function~$V_N$ in view of the finite prediction horizon~$N$.
\begin{proposition}\label{prop:link_costcont}
	Let the error bound~\eqref{eq:operatorbound} hold with $\varepsilon > 0$ and the stage cost be given by~\eqref{eq:stage_cost}.
	Suppose existence of a monotonically increasing and bounded sequence $(B_k)_{k \in \mathbb{N}} \subset \mathbb R$ and a set $S \subseteq \bX \ominus \mathcal{B}_\varepsilon(0)$ such that the growth bound
	\begin{align}\label{eq:cost_controllability}
	V_k(\hat{x}) \leq J_k(\hat{x},\hat{u}) \leq B_k \ell^\star(\hat{x}) \qquad\forall\,k\in \mathbb{N}
	\end{align}
	with $\ell^\star(\hat{x}) : = \inf_{u \in \bU} \ell(\hat{x},u)$ holds for all $\hat{x} \in S$ and some $\hat{u} = \hat{u}(\hat{x}) \in \mathcal{U}_N(\hat{x}) \cap \mathcal{U}_N^\varepsilon(\hat{x})$.
	Then, there exists a monotonically increasing and bounded sequence $(B_k^\varepsilon)_{k \in \mathbb{N}} \subset \mathbb{R}$ such that Inequality~\eqref{eq:cost_controllability} holds for $V_k^\varepsilon$ and $J_k^\varepsilon$ instead of $V_k$ and $J_k$, respectively. Moreover, we have $B_k^\varepsilon \rightarrow B_k$ for $\varepsilon \rightarrow 0$, $k \in \mathbb{N}$.
\end{proposition}
\begin{proof}
	Let $\tilde{x}(\cdot)$ and ${x}(\cdot)$ denote the trajectories generated by $\tilde{x}(n+1) = f^\varepsilon(\tilde{x}(n),\hat{u}_n)$ and $x(n+1) = f(x(n),\hat{u}_n)$, $n \in \mathbb{N}_0$, with $\tilde{x}(0) = \hat{x} = x(0)$, respectively. Set $\bar{\lambda} = \max\{|\lambda| : \lambda\ \mathrm{eigenvalue\ of}\ R\ \mathrm{or}\ Q\}$ and $0<\underline{\lambda} = \min\{|\lambda| : \lambda\ \mathrm{eigenvalue\ of}\ R\ \mathrm{or}\ Q\}$.
	Then, we have
	\begin{align}
	 \ell(\tilde{x}(n),\hat{u}_n) &= \| (\tilde{x}(n) - x(n)) + x(n) \|^2_Q + \| \hat{u}_n \|^2_R \label{eq:cost_controllability:proof:1} \\
	&\leq  \bar{\lambda} \| \tilde{x}(n) - x(n) \|^2 + 2 \bar{\lambda} 
	\| \tilde{x}(n) - x(n) \| \| x(n) \|+ \ell(x(n),\hat{u}_n), \nonumber
	\end{align}    
	If~\eqref{eq:operatorbound} holds, then Proposition~\ref{prop:errbound} yields the bound \eqref{eq:dynamics_bound} on the difference of~$f$ and~$f^\varepsilon$. Thus, we may estimate the term $e_{n+1} := \| \tilde{x}(n+1) - x(n+1) \|$ by
	\begin{align*}
	e_{n+1} & = \| f^\varepsilon(\tilde{x}(n),\hat{u}_{n}) \pm f(\tilde{x}(n),\hat{u}_{n}) - f(x(n),\hat{u}_{n}) \| \\
	& \leq \varepsilon \left( L_\Psi \| \tilde{x}(n) \mp x(n) \| + \Delta t \tilde{c} \| \hat{u}_n \| \right) + L_f e_n \\ 
	& = \varepsilon \cdot \bar{c} \left(\| x(n) \| + \| \hat{u}_n \| \right) + d e_n
	\end{align*}
	with $\bar{c} := \max \{ L_\Psi, \Delta t \tilde{c} \}$ and $d := L_f + \varepsilon L_\Psi$. Hence, 
	\begin{align}
	e_{n}^2 & \leq 4 \varepsilon^2 \bar{c}^2 ( \| x(n-1) \|^2 + \| \hat{u}_{n-1} \|^2 ) + 2 d^2 e_{n-1}^2 \nonumber \\
	& \leq \frac {4 \varepsilon^2 \bar{c}^2}{\underline{\lambda}} \sum_{i=0}^{n-1} (2d^2)^{n-1-i} \ell(x(i),\hat{u}_{i}), \nonumber \\
	e_{n} \| x(n) \hspace*{-0.25mm} \| \hspace*{-0.5mm} & \leq \varepsilon \bar{c} \text{\footnotesize $\frac{\| x(n\hspace*{-0.5mm}-\hspace*{-0.5mm}1) \hspace*{-0.125mm} \|^2 + \| \hat{u}_{n\hspace*{-0.125mm}-\hspace*{-0.125mm}1} \hspace*{-0.125mm} \|^2 + 2\| x(n) \hspace*{-0.25mm} \|^2}{2}$} \hspace*{-0.25mm}+\hspace*{-0.25mm} d e_{n\hspace*{-0.125mm}-\hspace*{-0.125mm}1} \| x(n) \hspace*{-0.25mm} \| \nonumber \\
	& \leq \frac {\varepsilon \bar{c}}{2 \underline{\lambda}} \sum_{i=0}^{n-1} d^{n-1-i} \Big( \ell(x(i),\hat{u}_i) + \ell^\star(x(n)) \Big) \nonumber 
	\end{align}
	Summing up the resulting inequalities for $\ell(\tilde{x}(n),\hat{u}_n)$ over $n \in [1:N-1]$ and using that the first summands in $\tilde{J}_N(\hat{x},\hat{u})$ and $V_N(\hat{x})$ coincide, we get 
	\begin{align*}
	\tilde{V}_N(\hat{x}) & \leq \tilde{J}_N(\hat{x},\hat{u}) \stackrel {\eqref{eq:cost_controllability:proof:1}}{\leq} J_N(\hat{x},\hat{u}) + \bar{\lambda} \left( \sum_{n=1}^{N-1} e_n^2 + 2 e_n \| x(n) \| \right) \\
	& \leq \left( B_N 
	+ \varepsilon \frac{ \bar{c} \bar{\lambda} }{2\underline{\lambda}} c_1 + \varepsilon^2 \frac {4 \bar{c}^2 \bar{\lambda}}{\underline{\lambda}} c_2 \right) \ell^\star(\hat{x}) =: B_N^\varepsilon \ell^\star(\hat{x}) \nonumber 
	\end{align*}
	with constants $c_1 = \sum_{n=1}^{N-1} d^{n-1} B_n + d^{N-1} B_N$ and $c_2 = \sum_{n=1}^{N-1} (2d^2)^{n-1} B_n$, where we have invoked the imposed cost controllability multiple times.
\end{proof}

Finally, invoking our findings on cost controllability, we verify the remaining conditions of Proposition~\ref{thm:PAS} to show the main result.
\begin{theorem}[PAS of EDMD-based MPC]\label{thm:main} 
	Let the error bound~\eqref{eq:operatorbound}, $\varepsilon \in (0,\varepsilon_0]$, for some $\varepsilon_0 > 0$, Assumption~\ref{ass:invariance} and cost controllability of the dynamics~\eqref{eq:dynamics_sampled_data} and the stage cost~\eqref{eq:stage_cost}, i.e., Condition~\eqref{eq:cost_controllability}, hold. 
	Let the prediction horizon~$N$ be chosen such that $\alpha \in (0,1)$ holds with
	\begin{equation}\label{eq:alpha}
	\alpha = \alpha_N := 1 - \frac {(B_2 - \omega) (B_N - 1) \prod_{i=3}^N (B_i - 1)} {\prod_{i=2}^N B_i - (B_2 - \omega) \prod_{i=3}^N (B_i - 1)}
	\end{equation}
	and $\omega = 1$.\footnote{The performance index or degree of suboptimality~$\alpha_N$ was proposed in~\cite{Grun09} and~\cite[Theorem~5.4]{GrunPann10} and updated to~\eqref{eq:alpha} in~\cite{Wort11}.}
	Further, let $S \subset \bX \ominus \mathcal{B}_{\varepsilon_0}(0)$ contain the origin in its interior and $\eta > 0$ be chosen such that, for all $\hat{x} \in S$, an optimal control function~$u^\star \in \mathcal{U}_N^\varepsilon(\hat{x})$ exists satisfying $x^\varepsilon_{u^\star}(k;\hat{x}) \in \bX \ominus \mathcal{B}_{\varepsilon + \eta}(0)$, $k\in [0:N-1]$.
	Then the EDMD-based MPC controller ensures semi-global practical asymptotic stability of the origin w.r.t.\ $\varepsilon$ on the set~$S$.
\end{theorem}
\begin{proof}
	First, we show condition~(i) of Proposition~\ref{thm:PAS} for the system dynamics~\eqref{eq:dynamics_approx}.  
	To this end, note that the lower bound on the optimal value function can be inferred by 
	\begin{align*}
	V_N^\varepsilon(\hat{x}) = \!\!\!\inf_{u\in \mathcal{U}_N^\varepsilon(\hat{x})} \!\!J_N^\varepsilon(\hat{x},u) \geq \!\inf_{u\in \mathbb{U}}\ell(\hat{x},u) = \| \hat{x} \|_Q^2 \geq \underline{\lambda} \| \hat{x} \|^2
	\end{align*}
	with $\underline{\lambda} > 0$ defined as in the proof of Proposition~\ref{prop:link_costcont}. 
	Then, defining $\alpha^{\varepsilon_0}$ analogously~\eqref{eq:alpha} using the sequence $(B_n^{\varepsilon_0})_{n = 2}^{N}$ instead and invoking $\lim_{\varepsilon_0 \searrow 0} B_n^{\varepsilon_0} = B_n$ yields $\alpha^{\varepsilon_0} \in (\alpha,1)$ for sufficiently small~$\varepsilon_0$. 
	This ensures the relaxed Lyapunov inequality for all $V_N^\varepsilon$, $\varepsilon \in (0,\varepsilon^0]$ by applying \cite[Theorem 5.2]{Grun09}. 
	Further, the upper bound on the value function $V_N^\varepsilon(\hat{x})$ directly follows from the imposed (and preserved) cost controllability. 
	Hence, we established the value function~$V_N^\varepsilon$ as a Lyapunov function for the closed loop of the surrogate dynamics~$f^\varepsilon$. 
	
	It remains to show $|V^\varepsilon_N(y_1)-V^\varepsilon_N(y_2)| \leq L\|y_1-y_2\|$ for all $y_1,y_2 \in S$ and $\varepsilon \in (0,\varepsilon_0]$, i.e., uniform continuity of~$V^\varepsilon_N$ with $\omega_V(r)=Lr$, for some $L > 0$.
	Then, the condition of Proposition~\ref{thm:PAS} hold and the assertion follows. 
	
	In combination with the uniform continuity of~$f^\varepsilon$ proven in Lemma~\ref{lem:as12}, the assumption $x^\varepsilon_{u^\star}(k;\hat{x}) \in \bX \ominus \mathcal{B}_{\varepsilon + \eta}(0)$ for all $k\in [0:N-1]$ implies the existence of~$\hat{\eta}>0$ such that, for each $\hat{x} \in S$, the respective optimal control $u^\star \in \mathcal{U}_N^\varepsilon(\hat{x})$ remains admissible for all initial values from $\mathcal{B}_{\hat{\eta}}(\hat{x})$. 
	Then, $V_N^\varepsilon$ is uniformly bounded on~$S$. This immediately shows the assertion for $y_1, y_2 \in S$ with $\| y_1 - y_2 \| > \hat{\eta}$, see, e.g., \cite{BoccGrun14} for a detailed outline of the construction. 
	Hence, it remains to show the assumption for $y_1,y_2 \in S$ satisfying $\| y_1 - y_2 \| \leq \hat{\eta}$. 
	Based on our assumption that an optimal sequence of control values exists, for every $y_2 \in \bX$ there is $u_2^{\star} \in \calU_N^\varepsilon(y_2)$ such that $V^\varepsilon_N(y_2)=J_N(y_2,u_2^{\star})$. 
	Then, invoking admissibility of $u_2^\star$ for~$y_1$, uniform Lipschitz continuity of $f^\varepsilon(\cdot,u)$ on $S$ in $\varepsilon \in (0,\varepsilon_0]$ and $u \in \mathbb{U}$, we get
	\begin{align*}
		 V^\varepsilon_N(y_1)-V^\varepsilon_N(y_2) &\leq J_N^\varepsilon(y_1,u_2^{\star}) - J_N^\varepsilon(y_2,u_2^{\star}) \\
		& = \sum_{k=0}^{N-1} \| x_{u_2^{\star}}^\varepsilon(k;y_1) - x_{u_2^{\star}}^\varepsilon(k;y_2) \|_Q^2 
	 + 2 x_{u_2^{\star}}^\varepsilon(k;y_2)^\top Q (x_{u_2^{\star}}^\varepsilon(k;y_1)-x_{u_2^{\star}}^\varepsilon(k;y_2)) \\
		& \leq  \bar{\lambda} \bar{c} \bigg[ \bar{c} \hat{\eta} + 2N \| x_{u_2^{\star}}^\varepsilon(k;y_2) \| \bigg) \bigg] \| y_1 - y_2 \|
	\end{align*}
	with $\bar{c} := \sum_{k=0}^{N-1} (c(\varepsilon_0) {L}_\Psi)^k$ for all $y_1,y_2\in \bX$. 
	Then, using that $\| x_{u_2^{\star}}^\varepsilon(k;y_2) \|$ is uniformly bounded on the compact set~$\bX$, we have derived $V^\varepsilon_N(y_1)-V^\varepsilon_N(y_2) \leq L \| y_1 - y_2 \|$.
	Analogously, 
	\begin{align*}
	V^\varepsilon_N(y_2) - V^\varepsilon_N(y_1) \leq J_N(y_2,u_1^{\star}) - J_N(y_1,u_1^{\star}) \leq L \|y_1-y_2\|
	\end{align*}
	on $S$. 
	Combining both inequalities yields the assertion.
\end{proof}

The assumption that the minimum exists may be completely dropped and is only imposed to streamline the presentation, see, e.g., \cite[p.\ 59]{GrunPann17} for details.
The imposed (technical) condition w.r.t.\ $\eta > 0$ can, e.g., be ensured by choosing a sufficiently small sub-level set $\{ x \in S : V_N^\varepsilon(x) \leq a\}$ such that $x_{u^\star}^\varepsilon(k) \notin \bX \ominus \mathcal{B}_{\varepsilon + \eta}(0)$ for some $k \in [1:N-1]$ yields a contradiction in view of the quadratic penalization of that state in the stage cost and the assumed bound~$a$ on the sub-level set --~similar to the construction used in~\cite{BoccGrun14}.

The assumed bound~\eqref{eq:operatorbound} of Theorem~\ref{thm:main} and cost controllability of the original system are the key ingredients for PAS of EDMD-based MPC. 
In Proposition~\ref{prop:generatorbound} we proved that such a bound can be guaranteed with probability $1-\delta$. This allows to also deduce PAS with probability $1-\delta$. Increasing the number of samples can then be used to either increase the confidence (that is, to reduce $\delta$), or reduce $\varepsilon$. The latter allows to shrink the set of PAS, i.e., reduce the radius $r>0$ in Definition~\ref{def:pracstab}.

\section{Numerical simulations}
\label{sec:example}

\noindent In this section we conduct numerical simulations to validate practical asymptotic stability of the origin for EDMD-based MPC as rigorously shown in Theorem~\ref{thm:main}. 

First, we consider the van-der-Pol oscillator given by
\begin{align} \label{eq:vanderPol}
\begin{pmatrix}
\dot x_1(t) \\ \dot x_2(t)
\end{pmatrix} = \begin{pmatrix}
x_2(t) \\ \mu (1 - x_1^2(t)) x_2(t) - x_1(t) + u(t)
\end{pmatrix}
\end{align}
for $\mu=0.1$. 
Since the linearization at the origin is controllable, cost controllablility holds for the quadratic stage cost~\eqref{eq:stage_cost}, see, e.g., \cite{Wort11}.
We consider the ODE~\eqref{eq:vanderPol} as a sampled-data system with zero-order hold as introduced in~\eqref{eq:dynamics_sampled_data}, where the integrals are numerically solved using the Runge-Kutta-Fehlberg method (RK45) with step-size control (Python function \textit{scipy.integrate.solve\_ivp}). 
For the approximation of the Koopman operator on the set $\mathbb{X} = [-2, 2]^2$, EDMD as described in Section~\ref{sec:eDMD} is used. As dictionary of observables we choose all $n_x$-variate monomials of degree less or equal than three, resulting in a dictionary size of $M = 10$. 
The step size is set to $\Delta t = 0.05$.

\begin{figure}[htb]
	\centering
	\includegraphics[width=.50\linewidth]{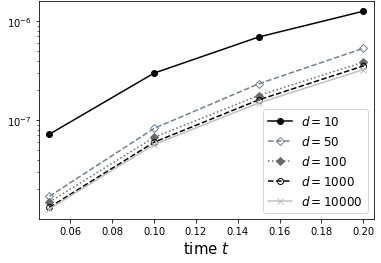} 
	\caption{Averaged error of the EDMD-based solution for different number of data points for fixed random control sequence.}
	\label{fig:openloop}
\end{figure}

First, we inspect the open-loop error of the EDMD-based surrogate for a random but fixed control sequence~$u$ and different numbers of data points $d \in \lbrace 10, 50, 100, 1000, 10000 \rbrace$, cf.\ Figure~\ref{fig:openloop}, which shows the average norm of the error for $100$ initial conditions distributed uniformly over the set $\mathbb{X} = [-2, 2]^2$. As to be expected from Proposition~\ref{prop:errbound}, the open-loop error decreases for increased number of samples. 

Next, we inspect the MPC closed-loop while imposing the constraints $-5 \leq u(k) \leq 5$ for $k \in [0:N-1]$ and $x(k)\in \bX$ for $k \in [0:N]$, respectively. 
We compare the closed-loop performance resulting from nominal MPC denoted by~$x_{\mu_N}$ as defined in~\eqref{eq:dynamics_closed_loop}
and EDMD-based MPC $x_{\mu_N^\varepsilon}$ defined in \eqref{eq:ex_cl} for $\lambda \in \lbrace 0.05, 0.25\rbrace$ and optimization horizons $N \in \lbrace 30, 50\rbrace$. 
The Koopman approximation is performed using EDMD with $d = 10000$ i.i.d.\ data points.
For small control penalization parameter $\lambda = 0.05$, the norm of the closed-loop state corresponding to nominal MPC decays until the precision $10^{-12}$ of the optimization solver is reached. As to be expected, this decay is faster for a longer prediction horizon. 
As proven in Theorem~\ref{thm:main}, the EDMD-based surrogate only enjoys practical asymptotic stability. 
More precisely, increasing the horizon only increases the convergence speed, but does not lead to a lower norm at the end of the considered simulation horizon. 

\begin{figure}[htb]
	\centering
	\includegraphics[width=.50\linewidth]{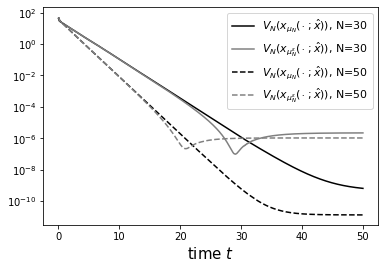}
	\caption{Optimal value functions along the closed-loop of system~\eqref{eq:vanderPol} for nominal MPC~(black) and EDMD-based MPC (gray) for horizons $N = 30$ (solid) and $N = 50$ (dashed) for $\lambda = 0.25$.}
	\label{fig:VN vanderpol}
\end{figure}

In Figure~\ref{fig:VN vanderpol}, we illustrate the decrease of the optimal value function along the closed-loop trajectories.
The observed stagnation indicates that the bottleneck is the approximation quality of the EDMD-based surrogate. 
The behavior is qualitatively very similar to the norm of the solution. 
Moreover, we observe a strict decrease of the value function over time. This is not the case for the EDMD-based MPC, for which we only have \emph{practical} asymptotic stability of the origin. 
Correspondingly, $V_N(x_{\mu_N^{\varepsilon}}(\cdot;\hat{x}))$ only decreases outside of a neighboorhood of the origin. 

\noindent 
The next example is taken from~\cite{NaraKwon19}, where the parameter values can be found. Here, $\dot{x}(t) = f(x(t), Q)$ describes an exothermic reaction that converts reactant $A$ to product $B$ and is given by
\begin{align}\label{eq:chemical}
\dot{x}(t) &= \begin{pmatrix}
\frac{F}{V_r}(C_{A0} - C_A) - k_0e^{\frac{-E}{RT_r}}C_A^2 \\
\frac{F}{V_r}(T_{A0} - T_r) -\frac{\Delta H}{\rho C_p}k_0e^{\frac{-E}{RT_r}}C_A^2+\frac{Q}{\rho C_p V_r}
\end{pmatrix}
\end{align}
with state~$x = (C_A, T_r) ^\top\in \R^2$, where $C_A$ is the concentration of $A$, $T_r$ the reactor temperature, and the control input~$Q$ is the heat supplied to the reactor. 
Since we want to stabilize the controlled steady state~$x^s = (C_{As}, T_{rs})^\top = (1.907, 300.6287)^\top$ ($Q^s = 0 \ \mathrm{kJ/hr}$), we consider the shifted dynamics, for which is origin is a steady state.

For EDMD, we use $d = 1000$ i.i.d.\ data points~$x^i$ drawn from the state-constrained set $\mathbb{X} = [-0.5, 0.5] \times [-20, 30]$ and propagate them by $\Delta t = 10^{-2}$ time units for control input $u_0 = 0$ and $u_1 = 1000$, respectively. 
The dictionary consists of observables $\{1, x_1, x_2, x_1^2, x_2^2, e^{\nicefrac{1}{x_1}}, e^{\nicefrac{1}{x_2}}\}$. 
We consider the respective OCP subject to $\mathbb{U} = [-10000,10000]$ with weighting $\lambda = 10^{-6}$ and $P = \mathrm{diag}(10^2, 1)$ for control and state.

\begin{figure}
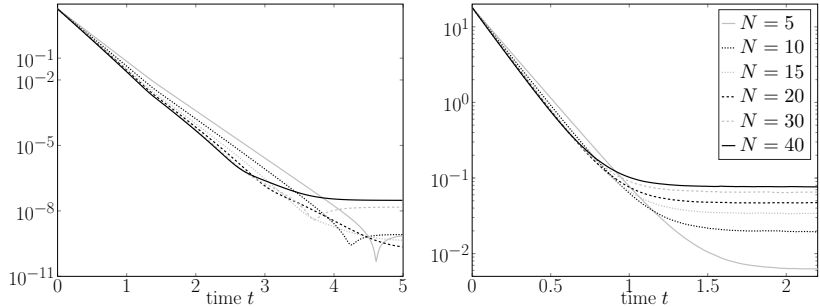

	\centering
	\input{6horizons_for_ODE}
	\input{5horizons_for_d1000}
	\caption{MPC closed loop $\| x_{\mu_N}(\cdot;x_0) \|$ (left) and $\| x_{\mu_N}^\varepsilon(\cdot;x_0) \|$ (right, EDMD: $d = 1000$) for system dynamics~\eqref{eq:chemical} for different horizons~$N$.} 
	\label{fig:CSTR}
\end{figure}

Figure~\ref{fig:CSTR} shows the numerical simulations emanating from the initial condition $x_0 = (0.5, -18)^\top$. The decay in norm of the closed-loop state corresponding to the EDMD-based surrogate stagnates around $10^{-2}$, i.e., practical asymptotic stability can be observed in this example, too. For the considered horizons, the decreasing behavior in the beginning until the point of stagnation is reached is similar to that of nominal MPC. The fact that the convergence stagnates earlier for larger $N$ is not unexpected, because $\omega_V$ in Proposition \ref{thm:PAS}(ii) may deteriorate since larger $N$ may render the optimal values more sensitive w.r.t.\ the initial condition.

\section{Conclusions}
\label{sec:conclusions}

\noindent 
We proved practical asymptotic stability of data-driven MPC for nonlinear systems using EDMD embedded in the Koopman framework. 
To this end, we established a novel bound on the estimation error, which scales proportional to the norm of the state \textit{and} the control. The underlying idea of imposing a certain structure in EDMD and, then, deriving \textit{proportional} bounds was also key in follow-up work for controller design using the Koopman generator~\cite{StraScha23} and operator~\cite{StraScha24}. 
Then, we showed that cost controllability of the original model is preserved for the proposed data-based surrogate. Last, we provided two numerical examples to illustrate our findings and, in particular, the practical asymptotic stability of the origin.
	\bibliographystyle{abbrv}
	\bibliography{references}

\begin{thebibliography}{10}

\bibitem{BoccGrun14}
A.~Boccia, L.~Gr{\"u}ne, and K.~Worthmann.
\newblock Stability and feasibility of state constrained {MPC} without
  stabilizing terminal constraints.
\newblock {\em Systems \& Control letters}, 72:14--21, 2014.

\bibitem{BrudFu21}
D.~Bruder, X.~Fu, and R.~Vasudevan.
\newblock Advantages of bilinear {K}oopman realizations for the modeling and
  control of systems with unknown dynamics.
\newblock {\em IEEE Robotics Automat. Lett.}, 6(3):4369--4376, 2021.

\bibitem{BrunBrun16}
S.~L. Brunton, B.~W. Brunton, J.~L. Proctor, and J.~N. Kutz.
\newblock {K}oopman invariant subspaces and finite linear representations of
  nonlinear dynamical systems for control.
\newblock {\em PloS one}, 11(2):e0150171, 2016.

\bibitem{CoroGrun20}
J.-M. Coron, L.~Grüne, and K.~Worthmann.
\newblock Model predictive control, cost controllability, and homogeneity.
\newblock {\em SIAM Journal on Control and Optimization}, 58(5):2979--2996,
  2020.

\bibitem{EsteWort20}
W.~Esterhuizen, K.~Worthmann, and S.~Streif.
\newblock Recursive feasibility of continuous-time model predictive control
  without stabilising constraints.
\newblock {\em IEEE Control Systems Letters}, 5(1):265--270, 2020.

\bibitem{FolkBurd21}
C.~Folkestad and J.~W. Burdick.
\newblock {K}oopman {NMPC}: {K}oopman-based learning and nonlinear model
  predictive control of control-affine systems.
\newblock In {\em IEEE International Conference on Robotics and Automation
  (ICRA)}, pages 7350--7356, 2021.

\bibitem{GoswPale21}
D.~Goswami and D.~A. Paley.
\newblock Bilinearization, reachability, and optimal control of control-affine
  nonlinear systems: A {K}oopman spectral approach.
\newblock {\em IEEE Trans. Automat. Control}, 67(6):2715--2728, 2021.

\bibitem{Grun09}
L.~Gr{\"u}ne.
\newblock Analysis and design of unconstrained nonlinear {MPC} schemes for
  finite and infinite dimensional systems.
\newblock {\em SIAM Journal on Control and Optimization}, 48(2):1206--1228,
  2009.

\bibitem{GrunPann17}
L.~Gr{\"u}ne and J.~Pannek.
\newblock {\em Nonlinear model predictive control}.
\newblock Springer Cham, 2017.

\bibitem{GrunPann10}
L.~Gr{\"u}ne, J.~Pannek, M.~Seehafer, and K.~Worthmann.
\newblock Analysis of unconstrained nonlinear {MPC} schemes with time varying
  control horizon.
\newblock {\em SIAM Journal on Control and Optimization}, 48(8):4938--4962,
  2010.

\bibitem{HaseCort21}
M.~Haseli and J.~Cort{\'e}s.
\newblock Learning {K}oopman eigenfunctions and invariant subspaces from data:
  {S}ymmetric subspace decomposition.
\newblock {\em IEEE Transactions on Automatic Control}, 67(7):3442--3457, 2021.

\bibitem{KanaYama22}
M.~Kanai and M.~Yamakita.
\newblock Linear model predictive control with lifted bilinear models by
  {K}oopman-based approach.
\newblock {\em SICE Journal of Control, Measurement, and System Integration},
  15(2):162--171, 2022.

\bibitem{KordMezi18b}
M.~Korda and I.~Mezi{\'c}.
\newblock Linear predictors for nonlinear dynamical systems: {K}oopman operator
  meets model predictive control.
\newblock {\em Automatica}, 93:149--160, 2018.

\bibitem{KordMezi18a}
M.~Korda and I.~Mezi{\'c}.
\newblock On convergence of extended dynamic mode decomposition to the
  {K}oopman operator.
\newblock {\em Journal of Nonlinear Science}, 28(2):687--710, 2018.

\bibitem{KordMezi20}
M.~Korda and I.~Mezi{\'c}.
\newblock Optimal construction of {K}oopman eigenfunctions for prediction and
  control.
\newblock {\em IEEE Transactions on Automatic Control}, 65(12):5114--5129,
  2020.

\bibitem{KrolTell22}
A.~Krolicki, D.~Tellez-Castro, and U.~Vaidya.
\newblock Nonlinear dual-mode model predictive control using {K}oopman
  eigenfunctions.
\newblock In {\em 61st IEEE Conference on Decision and Control (CDC)}, pages
  3074--3079, 2022.

\bibitem{LasoMack13}
A.~Lasota and M.~C. Mackey.
\newblock {\em Chaos, fractals, and noise: stochastic aspects of dynamics}.
\newblock Springer New York, 2013.

\bibitem{MaHuan19}
X.~Ma, B.~Huang, and U.~Vaidya.
\newblock Optimal quadratic regulation of nonlinear system using {K}oopman
  operator.
\newblock In {\em Proceedings of the 2019 IEEE American Control Conference
  (ACC)}, pages 4911--4916, 2019.

\bibitem{MamaCast21}
G.~Mamakoukas, M.~L. Castano, X.~Tan, and T.~D. Murphey.
\newblock Derivative-based {K}oopman operators for real-time control of robotic
  systems.
\newblock {\em IEEE Transactions on Robotics}, 37(6):2173--2192, 2021.

\bibitem{MaurSusu20}
A.~Mauroy, Y.~Susuki, and I.~Mezi{\'c}.
\newblock {\em {K}oopman operator in systems and control}.
\newblock Springer Cham, 2020.

\bibitem{NaraKwon19}
A.~Narasingam and J.~S.-I. Kwon.
\newblock Koopman lyapunov-based model predictive control of nonlinear chemical
  process systems.
\newblock {\em AIChE Journal}, 65(11):e16743, 2019.

\bibitem{NaraSon23}
A.~Narasingam, S.~H. Son, and J.~S.-I. Kwon.
\newblock Data-driven feedback stabilisation of nonlinear systems:
  {K}oopman-based model predictive control.
\newblock {\em International Journal of Control}, 96(3):770--781, 2023.

\bibitem{NuskPeit23}
F.~N{\"u}ske, S.~Peitz, F.~Philipp, M.~Schaller, and K.~Worthmann.
\newblock Finite-data error bounds for {K}oopman-based prediction and control.
\newblock {\em Journal of Nonlinear Science}, 33:14, 2023.

\bibitem{OttoRowl21}
S.~E. Otto and C.~W. Rowley.
\newblock {K}oopman operators for estimation and control of dynamical systems.
\newblock {\em Annual Review of Control, Robotics, and Autonomous Systems},
  4:59--87, 2021.

\bibitem{PeitOtto20}
S.~Peitz, S.~E. Otto, and C.~W. Rowley.
\newblock Data-driven model predictive control using interpolated {K}oopman
  generators.
\newblock {\em SIAM Journal on Applied Dynamical Systems}, 19(3):2162--2193,
  2020.

\bibitem{PhilScha2024}
F.~M. Philipp, M.~Schaller, K.~Worthmann, S.~Peitz, and F.~N{\"u}ske.
\newblock Error analysis of kernel {EDMD} for prediction and control in the
  {K}oopman framework.
\newblock {\em Preprint arXiv:2312.10460}, 2024.

\bibitem{PhilScha23}
F.~M. Philipp, M.~Schaller, K.~Worthmann, S.~Peitz, and F.~N{\"u}ske.
\newblock Error bounds for kernel-based approximations of the {K}oopman
  operator.
\newblock {\em Applied and Computational Harmonic Analysis}, 71:101657, 2024.

\bibitem{ProcBrun16}
J.~L. Proctor, S.~L. Brunton, and J.~N. Kutz.
\newblock Dynamic mode decomposition with control.
\newblock {\em SIAM J. Appl. Dynam. Syst.}, 15(1):142--161, 2016.

\bibitem{ProcBrun18}
J.~L. Proctor, S.~L. Brunton, and J.~N. Kutz.
\newblock Generalizing {K}oopman theory to allow for inputs and control.
\newblock {\em SIAM Journal on Applied Dynamical Systems}, 17(1):909--930,
  2018.

\bibitem{SchaWort23}
M.~Schaller, K.~Worthmann, F.~Philipp, S.~Peitz, and F.~N{\"u}ske.
\newblock Towards reliable data-based optimal and predictive control using
  extended {DMD}.
\newblock {\em IFAC-PapersOnLine}, 56(1):169--174, 2023.

\bibitem{SonNara22}
S.~H. Son, A.~Narasingam, and J.~S.-I. Kwon.
\newblock Development of offset-free {K}oopman {L}yapunov-based model
  predictive control and mathematical analysis for zero steady-state offset
  condition considering influence of {L}yapunov constraints on equilibrium
  point.
\newblock {\em Journal of Process Control}, 118:26--36, 2022.

\bibitem{StraBerb23}
R.~Str{\"a}sser, J.~Berberich, and F.~Allg{\"o}wer.
\newblock Control of bilinear systems using gain-scheduling: Stability and
  performance guarantees.
\newblock In {\em 62nd IEEE Conference on Decision and Control (CDC)}, pages
  4674--4681, 2023.

\bibitem{StraScha23}
R.~Str{\"a}sser, M.~Schaller, K.~Worthmann, J.~Berberich, and F.~Allg{\"o}wer.
\newblock Koopman-based feedback design with stability guarantees.
\newblock {\em IEEE Transactions on Automatic Control}, 2024.

\bibitem{StraScha24}
R.~Str{\"a}sser, M.~Schaller, K.~Worthmann, J.~Berberich, and F.~Allg{\"o}wer.
\newblock Saf{EDMD}: A certified learning architecture tailored to data-driven
  control of nonlinear dynamical systems.
\newblock {\em Preprint arXiv:2402.03145}, 2024.

\bibitem{Sura16}
A.~Surana.
\newblock {K}oopman operator based observer synthesis for control-affine
  nonlinear systems.
\newblock In {\em 55th IEEE Conference on Decision and Control (CDC)}, pages
  6492--6499, 2016.

\bibitem{WillHema2016}
M.~O. Williams, M.~S. Hemati, S.~T. Dawson, I.~G. Kevrekidis, and C.~W. Rowley.
\newblock Extending data-driven {K}oopman analysis to actuated systems.
\newblock {\em IFAC-PapersOnLine}, 49(18):704--709, 2016.

\bibitem{WillKevr15}
M.~O. Williams, I.~G. Kevrekidis, and C.~W. Rowley.
\newblock A data--driven approximation of the {K}oopman operator: Extending
  dynamic mode decomposition.
\newblock {\em Journal of Nonlinear Science}, 25:1307--1346, 2015.

\bibitem{Wort11}
K.~Worthmann.
\newblock {\em Stability analysis of unconstrained receding horizon control
  schemes}.
\newblock PhD thesis, University of Bayreuth, 2011.
\newblock
  \url{https://epub.uni-bayreuth.de/id/eprint/273/1/Dissertation_KarlWorthmann.pdf}.

\bibitem{YuShen22}
S.~Yu, C.~Shen, and T.~Ersal.
\newblock Autonomous driving using linear model predictive control with a
  {K}oopman operator based bilinear vehicle model.
\newblock {\em IFAC-PapersOnLine}, 55(24):254--259, 2022.

\bibitem{ZhanPan22}
X.~Zhang, W.~Pan, R.~Scattolini, S.~Yu, and X.~Xu.
\newblock Robust tube-based model predictive control with {K}oopman operators.
\newblock {\em Automatica}, 137:110114, 2022.

\end{thebibliography}

\end{document}